\documentclass[hein,final,leqno]{shein}
\usepackage{latexsym}
\usepackage{graphics}
\usepackage{epsfig}
\usepackage{epsf}
\usepackage{eucal}
\usepackage{pstcol}
\usepackage{amsfonts}
\usepackage{amstext}
\usepackage{amssymb}
\usepackage{amsmath}
\usepackage{amsbsy}
\usepackage{amsxtra}
\usepackage{amscd}
\newtheorem{ltheorem}{Theorem}           

\renewcommand{\abstractname}{Abstract}
\begin{document}
\title{DSC numerical solution of the Oberbeck-Boussinesq equations}
\author{Steffen Hein
}                     
%
%
\institute{SPINNER GmbH., M\"unchen, Theoretical Numerics [TB01],
Aiblinger Str. 30,\linebreak DE-83620 Westerham,
Germany  \hfill \email{s.hein@spinner.de}}
%
%
\date{}
%

\maketitle 
%
\begin{abstract}
Dual Scattering Channel schemes generalise \textsc{Johns}' TLM algorithm
and replace the latter in situations where the transmission line picture
of wave propagation fails. This is notoriously the case in applications
to fluid dynamics, for instance.
In this paper, a DSC numerical solution of the
\textsc{Oberbeck-Boussinesq} equations is presented,
which approximate the \textsc{Navier-Stokes} equations
for viscous quasi incompressible flow with moderate
variation in temperature.
\vspace{5pt} \\
\textbf{Keywords}:
Time domain methods, DSC schemes, fluid dynamics, CFD,
\textsc{Navier-Stokes} equations, \textsc{Boussinesq}
approximation.\\ 
{\vspace{.05cm}\hfill
\textbf{MSC-classes}:\;\textnormal{65C20,\,65M06,\,76D05}}

\vspace{.2cm}
\hspace{-.52cm}
\!\!Westerham on \today

\end{abstract}

\markboth{{\normalsize \textsc{Steffen Hein}}}
{{\normalsize \textsc{DSC numerical solution of Oberbeck-Boussinesq equations}}}

\normalsize
\vspace{-0.2cm}
\section{Introduction}
Crushing the continuum down into mesh cells
is a queer, artificial exercise.
There are yet natural ways of computing the fields
in a cellular mesh, and so to mitigate the desaster. 
Imposing a cellular mesh is tantamout to locally enforcing
cell-boundary duality upon space -~and the DSC setup
offers a natural way to handle such situations.

Dual Scattering Channel (DSC) schemes are characterized by a two-step
cycle of iteration which alternately updates the computed fields within
cells and on their interfaces. If the updating instructions are explicit,
then a near-field interaction principle leads to the typical structure
of the DSC algorithm, viz. to its scattering process interpretation.
A pair of vectors which represent the same field within a cell
and on its surface essentially constitutes a \emph{scattering channel}.
An equivalent definition can be given in terms of a pair of
distributions that 'measure' the field within the cell and on one of
its faces. In the \emph{primal} DSC scheme which is the Transmission
Line Matrix (TLM) method along \textsc{Johns}' line \cite{JoB}
ports of transmission line links visualise these distributions
(finite integrals, in this case). 

DSC schemes and their relations to the TLM method [Ch, Tlm1-3, Re]
have been conceptually and technically scrutinised in-depth \cite{He1}.
They are unconditionally stable under quite general circumstances
-~made tangible with the notion of \emph{$\alpha$-passivity} \cite{He2},
and they are especially suitable for handling boundary conditions, 
non-orthogonal mesh, or also for replacing a staggered grid where
otherwise need of such is.

In section~\ref{S:sec1}, we first recapitulate some characteristic
features of DSC schemes, which in extenso are treated in \cite{He1},
before we deal in section~\ref{S:sec2} with the
\textsc{Oberbeck - Boussinesq} approximation
to the \textsc{Navier-Stokes} equations.
The DSC model outlined here should be considered a \emph{prototype\/},
first of all.
In fact, the \textsc{Boussinesq} equations for viscous fluid flow,
inspite of retaining the usually predominant non-linear advective
part of the \textsc{Navier-Stokes} momentum equations,
can only claim limited range of validity,
due to their well known simplifications.
The \textsc{Oberbeck-Boussinesq} approach is, however, prototypical
also in providing a basis for many turbulence models \cite{ATP}
that can be implemented following essentially the lines of this paper.

The treated implementation with unstructured hexahedral mesh,
outlined in section~\ref{S:sec3}, has recently been implemented
and coupled to \textsc{Spinner}'s Maxwell field solver,
thus allowing for computing conductive and convective heat transfer
simultaneously with the electric and magnetic heat sources in a
lossy Maxwell field.
Besides the underlying ideas and simplifications that enter the
\textsc{Oberbeck-Boussinesq} approximation and its DSC formulation, 
some numerical results are displayed to illustrate and validate the
approach.
\vspace{-0.2cm}
\section{Elements of DSC schemes}\label{S:sec1}
In this section, we resume some typical features of DSC schemes, 
referring the technically interested reader always to the systematic
exposition \cite{He1}.

Given a mesh cell, we think of a \emph{port} as a vector valued
distribution, associated to a cell face (with support, however,
\emph{not} necessarily confined to that face, cf. sect.\ref{S:sec3}), 
which assigns a state vector $\,z^{\,p}\,=\,(\,p\,,\,Z\,)\,$
to a physical field $Z\,$, the latter represented by a smooth
vector valued function in space-time. We also require that in the
given cell a \emph{nodal image} ${p\sptilde\/}$ of $p\/$ exists,
such that
\vspace{-0.10cm}
\begin{equation}\centering\label{1.1}
(\,p{\;\sptilde},\,Z\,)\;
=\;(\,p\,\circ\,\sigma\,,\,Z\,)\;
=\;(\,p\,,\,Z\,\circ\,\sigma^{-1}\,)\;,
\end{equation}
for every $\,Z\,$ (\,of class $\,C^{\,\infty}\,$, e.g.\,)\;,
where $\,\sigma\,$ denotes the spatial translation
${\sigma:\mathbb{R}^3\,\to\,\mathbb{R}^3}$ 
that shifts the geometrical node (centre of cell) onto the
(centre of) the respective face, cf. Fig1.
\begin{figure}[!h]\centering
\vspace{-2.3cm}
\setlength{\unitlength}{.7cm}
\begin{pspicture}(0.0,0.0)(6.0,4.6)\centering
\psset{xunit=0.7cm,yunit=0.7cm}
\psline[linewidth=0.1mm]{-}(0.0,1.0)(2.0,0.0)
\psline[linewidth=0.1mm]{-}(2.0,0.0)(6.0,1.0)
\psline[linewidth=0.1mm]{-}(6.0,1.0)(3.0,1.5)
\psline[linewidth=0.1mm]{-}(3.0,1.5)(0.0,1.0)
\psline[linewidth=0.1mm]{-}(0.0,1.0)(0.0,2.0)
\psline[linewidth=0.1mm]{-}(2.0,0.0)(2.0,1.0)
\psline[linewidth=0.1mm]{-}(6.0,1.0)(6.0,2.0)
\psline[linewidth=0.1mm]{-}(3.0,1.5)(3.0,2.5)
\psline[linewidth=0.2mm]{->}(1.0,0.5)(4.5,1.25)
\psline[linewidth=0.2mm]{->}(1.5,1.25)(4.0,0.5)
\psline[linewidth=0.2mm]{->}(1.0,2.5)(4.5,3.25)
\psline[linewidth=0.2mm]{->}(1.5,3.25)(4.0,2.5)
\psline[showpoints=true,linewidth=0.4mm]{->}(2.75,2.875)(2.75,0.825)
\psline[showpoints=true,linewidth=0.4mm]{-}(2.75,0.875)(2.75,0.825)
\rput(2.6,3.3){\scriptsize node}
\rput(4.9,3.4){\small $p\sptilde$}
\rput(2.5,1.8){\small $\sigma$}
\rput(0.0,0.5){\scriptsize cell face}
\rput(4.9,1.4){\small $p$}
\end{pspicture}
\caption{\textsl{Port on a cell face with nodal image.}
\hfill}\label{F:1}
\end{figure}
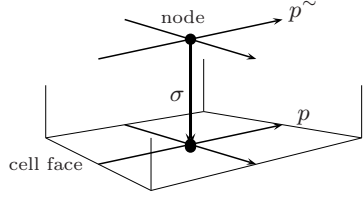

DSC fields split thus into \emph{port} and \emph{node} components,
$z^p$ and $z^n\,$, which represent the field at the cell surfaces
and within the cells.
The two components are updated at even and odd integer multiples,
respectively, of half a timestep ~${\tau}\,$ and are usually as step
functions constantly continued over the subsequent time intervals of
length $\tau\,$.
Moreover, we assume that the updating instructions are \emph{explicit\/},
i.e., with possibly time dependent functions ${F}$ and ${G}\,$, 
for $t=m\tau\,;\,m\in\mathbb{N}\;$
\begin{equation}\centering\label{1.2}
\begin{split}
\begin{aligned}
&z^n\,(\,t\,+\,\frac{\tau}{2}\,)\;
&&:\,=\quad F\,(\,[\,z^p\,]_{\,t\,}\,,
\;[\,z^n\,]_{\,t\,-\,\frac{\tau}{2}}\,)\;,\\
&z^p\,(\,t\,+\,\tau\,)\;
&&:\,=\quad G\,(\,[\,z^p\,]_{\,t}\,,
\;[\,z^n\,]_{\,t\,+\,\frac{\tau}{2}}\,)\;,
\end{aligned}
\end{split}
\end{equation}
where ${[\,z\,]_{\,t}}$ stands for the entire sequence up to time $\,t$
\begin{equation}\label{1.3}\notag
[\,z\,]_{\,t}\,:\,=\;(\,z\,(\,t\,-\,\mu\tau\,)\,)\,_{\mu\,\in\,\mathbb{N}}\quad
\end{equation}
and we agree upon fixing
${\,z^{\,p,n}\,(\,t\,)\,:\,=\,0\,}$ for ${\,t\,<\,0\,}$.
(\,The 'back in time running' form of the sequence has certain technical 
advantages, cf.~\cite{He1}.) 

A fundamental DSC principle is \emph{near-field interaction},
which spells that every updated state depends only on states
(up to present time $\,t\,$) in the immediate neighbourhood.
More precisely: The next nodal state depends only on states
(with their history) in the same cell and on its boundary,
and a subsequent port state depends only on states
(again with history) on the same face and in the adjacent nodes.

As a consequence of near-field interaction, every DSC process
allows for an interpretation as a multiple scattering process
in the following sense.

If $M$ is a mesh cell system and
${\,\partial\zeta}$
denotes the boundary of cell ${\zeta\in M}$, then every DSC state obviously
permits a unique \emph{scattering channel representation\/} in the space
\begin{equation}\centering\label{1.4}\notag
P\;:\,=\;\prod\nolimits_{\,\zeta\in M}\;
\prod\nolimits_{\,p\in\partial\zeta}\;
(\, z_{\,\zeta}^{\,p}\,,\;z_{\,\zeta}^{p\sptilde}\,)\;,
\end{equation}
with canonical projections
${\pi_{\,\zeta}^{\,p,\,n}\,:\,P\,\to\,P_{\,\zeta}^{\,p,\,n}}$
into port and node components of cell $\zeta\,$
(\,the cell index is omitted, in general,
if there is no danger of confusion\,).
Also, there is a natural involutary isomorphism 
${nb\,:\,P\,\to\,P\,}$
\begin{equation}\centering\label{1.5}\notag
nb\,:\,(\,z^{\,p}\,,\,z^{\,p\sptilde}\,)\,\mapsto\,
(\,z^{\,p\sptilde},\,z^{\,p}\,)\;,
\end{equation}
which is named the \emph{node-boundary map\/} and obviously 
maps ${\,P^{\,p}}$ onto ${\,P^{\,n}}$ and vice versa.
For every DSC process $\,z\,=\,{\,(\,z^{\,p}\,,\,z^{\,n}\,)(\,t\,)\,}$,
the following \emph{incident\/} and \emph{outgoing fields}
$\,{z_{\,in}^{\,p}}\,$ and $\,{z_{\,out}^{\,n}}\,$ 
are then recursively well (viz. uniquely) defined, and are processes in
${\,P^{\,p}\,}$ and ${\,P^{\,n}\,}$, respectively:\; \\ 
For $\,t < 0\,$, $\,{z_{\,in}^{\,p}\,(\,t\,)\,
:\,=\,z_{\,out}^{\,n}\,(\,t-\frac{\tau}{2}\,)\,:\,=\,0\;}\,$,
and for every $\,0\,\leq\,t\,=\,m\tau\,$;
$\,m\in\mathbb{N}\,$
\begin{equation}\centering\label{1.6}
\begin{split}
\begin{aligned}
z_{\,in}^{\,p}\,(\,t\,)\;
&:\,=\; z^{\,p}\,(\,t\,)\,
-\,nb\circ z_{\,out}^{\,n}\,(\,t\,-\,\frac{\tau}{2}\,)\;,\\
z_{\,out}^{\,n}\,(\,t\,+\,\frac{\tau}{2}\,)\;
&:\,=\; z^{\,n}\,(\,t+\,\frac{\tau}{2}\,)\,-\,nb\circ z_{\,in}^{\,p}\,(\,t\,)\;.
\end{aligned}
\end{split}
\end{equation}
Hence, at every instant holds $\,{z^{\,p}\,(\,t\,)}\,
=\,{nb\circ z_{\,out}^{\,n}\,(\,t\,-\,\frac{\tau}{2}\,)}\,
+\,{z_{\,in}^{\,p}\,(\,t\,)}\,$
and \\
$\,{z^{\,n}\,(\,t\,+\,\frac{\tau}{2}\,)}\,
=\,{nb\circ z_{\,in}^{\,p}\,(\,t\,)}\,
+\,{z_{\,out}^{\,n}\,(\,t\,+\,\frac{\tau}{2}\,)}\,$\,.
Then near-field interaction implies that every state is only a function of
states incident (up to present time $\,t\,$) on scattering channels
connected to the respective node or face.
\newline
Precisely, it is shown that
\begin{ltheorem}{\label{T:1}}.
A pair of functions
$\,\mathcal{R}\,$ and $\,\mathcal{C}\,$ exists, such that for every cell
$\,{\zeta\in M}\,$
the process
$\,z_{\,\zeta}^{\,n}\,=\,{\pi_{\,\zeta}^{\,n}\circ z}\,$
complies with
\vspace{-0.15cm}
\begin{equation}\centering\label{1.9}
z_{\,\zeta}^{\,n}\,(\,t\,+\,\frac{\tau}{2}\,)\;
=\;\mathcal{R}\,(\,(\,z_{\,in}^{\,p}\,(\,t\,-\,\mu\tau\,)\,
)_{\,p\in\partial\zeta\,;\;\mu\in\mathbb{N}}\,)\;
\end{equation}
\vspace{-0.15cm}
and the port process
$\,z_{\,\zeta}^{\,p}\,=\,\pi_{\,\zeta}^{\,p}\,\circ\,z\,$
satisfies
\begin{equation}\centering\label{1.10}
\begin{aligned}
z_{\,\zeta}^{\,p}\,(\,t\,+\,\tau\,)\;
=\;\mathcal{C}\,(\,(\,z_{\,out}^{\,n}\,(\,t\,+\,\frac{\tau}{2}-
\mu\tau\,)\,)_{\,n\,\mid\,\partial\zeta\,;\;\mu\in\mathbb{N}}\,)\;.\\
\scriptsize{(\;\text{'$\;\mid\,$' short-hand for 'adjacent to'}\;\;)}
\end{aligned}
\end{equation}
\end{ltheorem}
\vspace{-0.45cm}
\emph{Remarks}
\hspace{.65cm}

\begin{itemize}
\item[(i)]
The statements imply, of course, that
$\,z_{\zeta,\,out}^{\,n}\,$ and $\,z_{\zeta,\,in}^{\,p}\,$
are themselves functions of states incident
on connected scattering channels, since
\vspace{-.3cm}
\end{itemize}
\begin{equation}\centering\label{1.11}
\begin{aligned}
z_{\zeta,\,out}^{\,n}\,(\,t\,+\,\frac{\tau}{2}\,)\;
=\;\mathcal{R}\,(\,(\,z_{\,in}^{\,p}\,(\,t\,-\mu\tau\,)\,)_{\,p\,\in\,
\partial\zeta\,;\;\mu\in\mathbb{N}}\,)\,
-\,z_{\zeta,\,in}^{\,p}\,(\,t\,) \quad \text{and} \\
z_{\zeta,\,in}^{\,p}\,(\,t\,)\;
=\;\mathcal{C}\,(\,(\,z_{\,out}^{\,n}\,(\,t\,-\,\frac{\tau}{2}-\mu\tau\,)\,
)_{\,n\,\mid\,p\,;\;\mu\in\mathbb{N}}\,)\,
-\,z_{\zeta,\,out}^{\,p\sptilde}\,(\,t\,-\,\frac{\tau}{2}\,)
\end{aligned}
\end{equation}
\begin{itemize}
\item[(ii)]
$\,\mathcal{R}\,$ and $\,\mathcal{C}\,$ are named the
\emph{reflection} and \emph{connection\/} maps, respectively,
of the DSC algorithm.
\end{itemize}
\begin{itemize}
\item[(iii)]
Near field interaction implies computational stability, if the reflection
and connection maps are contractive or 
\nolinebreak{$\,\alpha$-\emph{passive}} \cite{He2},
in addition.
\end{itemize}

\begin{proof}
Mere translation of the \emph{near-field interaction} principle and induction:
The statement is trivial for $\,t\,<\,0\,$.
If, for instance, \eqref{1.9} holds until $\,t-\frac{\tau}{2}\,$,
then in virtue of \eqref{1.2} and remark (i) also
\begin{equation}\centering\label{1.12}\notag
\begin{split}
\begin{aligned}
&z_{\,\zeta}^n\,(\,t\,+\,\frac{\tau}{2}\,)\,
&&=\quad F\,\;(\;[\hspace{-47pt}\underbrace{\;z_{\zeta}^{\,p}\;}_{\hspace{-5pt}
\,z_{\zeta,\,in}^{\,p}\,(\,t\,)\,+\,nb\,z_{\zeta,\,out}^{\,n}\,(\,t\,
-\,\frac{\tau}{2}\,)\,,\qquad}\hspace{-47pt}]_{\,t}\,,\;[\;\hspace{-28pt}
\underbrace{z_{\zeta}^{\,n}}_{\hspace{+35pt}\,
\mathcal{R}\,[\,z_{\zeta,\,in}^{\,p}\,]_{\,t\,-\tau}}\,\hspace{-27pt}]_{
\,t-\frac{\tau}{2}}\,)
\end{aligned}
\end{split}
\end{equation}
is a function of states incident from $\,p\,\in\,\partial\zeta\,$
until $\,t\,$.
\end{proof}
\vspace{-0.2cm}
\section{The dynamic equations}\label{S:sec2}
DSC algorithms are thus simply characterized as two-step explicit
schemes that alternately update states in ports and nodes of a 
cellular mesh and which, in virtue of a near-field interaction principle,
allow for a canonical interpretation as scattering processes.
The latter exchange incident and reflected quantities between
cells and their interfaces.

The ports and nodes are related to physical fields by vector valued
distributions that evaluate the fields at cell faces and within the
cells of a cellular mesh.
Such a distribution may be a \emph{finite integral},
as in the case of the TLM method, where finite path integrals over
electric and magnetic fields are evaluated in a discrete approximation
to Maxwell's integral equations \cite{He3}.
In the simplest case, it is only a \emph{Dirac measure} that pointwise
evaluates a field (\,or a field component\,) within a cell and on its
surface. The distribution can also be a composite of Dirac measures that
evaluate a field at different points in the cell -~which applies,
for instance, to the gradient functional treated in sect.~\ref{S:sec3}.

Classical thermodynamics with, in particular, energy conservation
entail the \emph{convection-diffusion} equation for the temperature
$\,T\,$ in a fluid of velocity $\,\vec{u}\,$ with constant thermal
diffusivity $\,\alpha\,$, heat source(s) $\,q\,$,
and negligible viscous heat dissipation, viz.
\begin{equation}\centering\label{2.1}
\frac{\partial\,T}{\partial\,t}\;
+\;\vec{u}\,\cdot\,grad\;T\;
=\;\alpha\,\Delta\,T\;+\;q\;.
\end{equation}
This is the energy equation for \emph{Boussinesq-incompressible\/}
fluids, e.g. \cite{GDN}.
The Navier-Stokes \emph{momentum equations\/} for a fluid of
dynamic viscosity $\,\mu\,$,
under pressure $\,p\,$,
and in a gravitational field of acceleration $\,\vec{g}\,$ require
\begin{equation}\centering\label{2.2}
\frac{\partial}{\partial\,t}\,(\,\varrho\,\vec{u}\,)\;
+\;(\,\vec{u}\cdot grad\,)\,(\,\varrho\,\vec{u}\,)
+\;grad\,p\;
=\;\mu\,\Delta\,\vec{u}\;
+\;\varrho\;\vec{g}\;.
\end{equation}
The \emph{Oberbeck-Boussinesq\/} approximation \cite{Obb},\cite{Bss}
starts from the assumption that the fluid properties are constant,
except fluid density, which only in the gravitational term varies
linearly with temperature;
and that viscous dissipation can be neglected.
Equations \eqref{2.2} become so with
$\,\varrho_{\,\infty}=const\,$ and
$\,\varrho\,(\,T\,)=\varrho_{\,\infty}\,\beta\,
(\,T\,(\,t\,,\vec{x}\,)\,-\,T_{\,\infty}\,)\;$; $\,\beta\,
:\,=\,\varrho^{-1}\,\partial\,\varrho\,/\,\partial\,T\,$
\begin{equation}\centering\label{2.3}
\frac{\partial\,\vec{u}}{\partial\,t}\,
+\,(\vec{u}\cdot grad\,)\,\vec{u}\,
+\,\frac{grad\,p}{\varrho_{\infty}}\,
=\,\frac{\mu}{\varrho_{\infty}}\,\Delta\,\vec{u}\;
+\,\beta\,(\,T\,(\,t,\,\vec{x}\,)
- T_{\infty}\,)\,\vec{g}\,.
\end{equation}
With the Gauss-Ostrogradski theorem applied to the
integrals over $\,{\Delta}\,=\,{div\,grad\,}\,$
on cell $\zeta$ with boundary ${\partial\zeta}$,
equations (\ref{2.1}, \ref{2.3}) yield,
with a time increment $\tau\,$,
the following updating instructions for $T$ and $\vec{u}$
averaged over the cell volume ${V_{\zeta}}$
\begin{equation}\centering\label{2.4}
\begin{aligned}
&T\,(\,t+\frac{\tau}{2}\,)\;
:\,=\\
&=\;\;T\,+\,\tau\,\{\,-\vec{u}\,\cdot\,grad\,T\,
+\,\frac{\alpha}{V_{\zeta}}\,\int\nolimits_{\partial\,\zeta}grad\,T\cdot dF\,
+\,\frac{1}{V_{\zeta}}\,\int\nolimits_{\zeta}\,q\,d V\;\}
\end{aligned}
\end{equation}
and
\begin{equation}\centering\label{2.5}
\begin{aligned}
&\vec{u}\,(\,t+\frac{\tau}{2}\,)\;:\,=\,\vec{u}\,
-\,\tau\,\{\,(\,\vec{u}\,\cdot\,grad\,)\,\vec{u}\,
+\,\frac{grad\,p}{\varrho_{\infty}}\,\}\;\;+\\
&\qquad+\;\tau\,\{\,\frac{\mu}{V_{\zeta}\,\varrho_{\infty}}\,
\int\nolimits_{\partial\,\zeta} grad\,\vec{u}\cdot dF\,
+\,\beta\,(\,T\,-\,T_{\infty}\,)\,\vec{g}\,\}\;.
\end{aligned}
\end{equation}
At the right-hand sides enter, of course, the last former updates
(at time $\,t\,-\,\tau/2\,$ and $\,t\,$, respectively)
of the nodal and cell face quantities.

$\,T\,$ and $\,\vec{u}\,$ are so updated at the reflection
step of the DSC algorithm. In contrast, the cell surface integrals
at the right-hand sides, in particular the gradients that
enter these, are updated on the connection step.
\newline
The next section proceeds with that in unstructured hexahedral mesh.
\vspace{-0.2cm}
\section{The non-orthogonal hexahedral cell}\label{S:sec3}
The physical interpretation of a DSC algorithm associates a
smoothly varying, i.e. in time and space sufficiently often
continuously differentiable (\,for instance, $\,C^{\,\infty}\,$\,)
scalar or vector field $\,Z\,$ to port and node states
$\,z^{\,p}\,$ and $\,z^{\,n}\,$ of a mesh cell system.

For notational economy (\,so avoiding many '$\sum\,$'\,s\,)
in the following we adopt \textsc{Einstein}'s convention to sum up over
identical right-hand sub and superscripts within terms where such are present,
while summation is \emph{not} carried out whenever a sub or superscript
also appears somewhere as a left-hand index (\,for instance, in
${(-1)}^{\kappa}\,a_{\kappa}^{\lambda}\,b_{\lambda}\,\,_{\kappa}c\;$
the sum is made over $\,\lambda\,$ but not over $\,\kappa\,$).

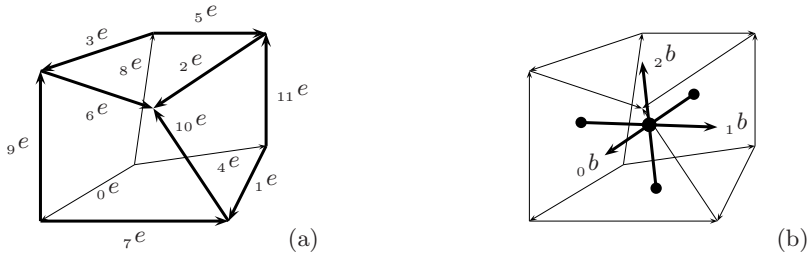
\begin{figure}[!h]\centering
\setlength{\unitlength}{1.cm}
\begin{pspicture}(-1.4,-.5)(20,3.5)\centering
\psset{xunit=.5cm,yunit=.5cm}
\psline[linewidth=0.1mm]{->}(2.5,1.5)(0.0,0.0)
\psline[linewidth=0.1mm]{->}(2.5,1.5)(6.0,2.0)
\psline[linewidth=0.1mm]{->}(2.5,1.5)(3.0,5.0)
\psline[linewidth=0.4mm]{->}(0.0,0.0)(0.0,4.0)
\psline[linewidth=0.4mm]{->}(0.0,0.0)(5.0,0.0)
\psline[linewidth=0.4mm]{->}(6.0,2.0)(6.0,5.0)
\psline[linewidth=0.4mm]{->}(6.0,2.0)(5.0,0.0)
\psline[linewidth=0.4mm]{->}(5.0,0.0)(3.0,3.0)
\psline[linewidth=0.4mm]{->}(3.0,5.0)(0.0,4.0)
\psline[linewidth=0.4mm]{->}(3.0,5.0)(6.0,5.0)
\psline[linewidth=0.4mm]{->}(0.0,4.0)(3.0,3.0)
\psline[linewidth=0.4mm]{->}(6.0,5.0)(3.0,3.0)
\rput(1.8,0.8){$_{_{0}} e$}
\rput(6.0,1.0){$_{_{1}} e$}
\rput(4.0,4.2){$_{_{2}} e$}
\rput(1.5,4.9){$_{_{3}} e$}
\rput(5.0,1.5){$_{_{4}} e$}
\rput(4.4,5.5){$_{_{5}} e$}
\rput(1.5,2.9){$_{_{6}} e$}
\rput(2.5,-.5){$_{_{7}} e$}
\rput(2.4,4.1){$_{_{8}} e$}
\rput(-.6,2.0){$_{_{9}} e$}
\rput(4.0,2.6){$_{_{10}} e$}
\rput(6.7,3.4){$_{_{11}} e$}
\rput(7.0,-0.5){\small{\textnormal{(a)}}}
\psline[linewidth=0.1mm]{->}(15.5,1.5)(13.0,0.0)
\psline[linewidth=0.1mm]{->}(15.5,1.5)(19.0,2.0)
\psline[linewidth=0.1mm]{->}(15.5,1.5)(16.0,5.0)
\psline[linewidth=0.1mm]{->}(13.0,0.0)(13.0,4.0)
\psline[linewidth=0.1mm]{->}(13.0,0.0)(18.0,0.0)
\psline[linewidth=0.1mm]{->}(19.0,2.0)(19.0,5.0)
\psline[linewidth=0.1mm]{->}(19.0,2.0)(18.0,0.0)
\psline[linewidth=0.1mm]{->}(18.0,0.0)(16.0,3.0)
\psline[linewidth=0.1mm]{->}(16.0,5.0)(13.0,4.0)
\psline[linewidth=0.1mm]{->}(16.0,5.0)(19.0,5.0)
\psline[linewidth=0.1mm]{->}(13.0,4.0)(16.0,3.0)
\psline[linewidth=0.1mm]{->}(19.0,5.0)(16.0,3.0)
\psline[showpoints=true,linewidth=0.4mm]{->}(17.375,3.375)(15.0,1.75)
\psline[showpoints=true,linewidth=0.4mm]{->}(14.375,2.625)(18.0,2.50)
\psline[showpoints=true,linewidth=0.4mm]{->}(16.37,.875)(16.0,4.25)
\psline[showpoints=true,
linewidth=0.6mm]{-}(16.1875,2.5625)(16.1875,2.5625) 
\rput(14.55,1.50){$_{_{0}} b$}
\rput(18.5,2.60){$_{_{1}} b$}
\rput(16.6,4.40){$_{_{2}} b$}
\rput(20,-0.5){\small{\textnormal{(b)}}}
\end{pspicture}
\caption{\textsl{Non-orthogonal hexahedral mesh cell. \newline
\textnormal{(a)} Edge vectors. \qquad
\textnormal{(b)} Node vectors.}\hfill}\label{F:2}
\end{figure}

Let a hexahedral cell be given by its eight vertices.
Define then \emph{edge vectors} ${(_{\nu} e)_{\nu=0,...,11}}\/$,
\emph{node vectors\/} ${(_{\mu} b)_{\mu = 0,1,2}}$,
and \emph{face vectors\/} ${(_{\iota} f)_{\iota = 0,...,5}}$,\,
using the labelling scheme of figure~\ref{F:2}\,a
\begin{equation}\label{3.1}\centering
\begin{split}
\begin{aligned}
_{\mu}b\;&:\,=\quad\frac{1}{4}
&& \!\!\sum\nolimits_{\nu = 0}^{3}\,_{_{(4\mu+\nu)}}e\,\,
&&\mu\,=\,0,1,2\\
\text{and}\quad
_{\iota}f\;&:\,=\;\,\frac{(-1)^{\,\iota}}{4}
&&\,(\,\,_{_{(8+2\iota)}}e\,
+_{_{(9+2(\iota+(-1)^{\iota}))}}e\,)
\,\,\land &&\\
& &&\;\;\land\,(\,_{_{(4+2\iota)}}e\, 
+_{_{(5+2\iota)}}e\,)\,\,&&\iota\,=\,0,...,5\,,
\end{aligned}
\end{split}
\end{equation}
with all indices understood cyclic modulo 12\, and $\,\land\,$
denoting the cross product in $\mathbb{R}^3$.

\begin{figure}[!h]\centering
\setlength{\unitlength}{1.0cm}
\begin{pspicture}(-2.4,-1.8)(10,4.0)\centering
\psset{xunit=.45cm,yunit=.45cm}
\psline[linewidth=0.2mm]{->}(0.0,4.0)(0.0,1.0)
\psline[linewidth=0.2mm]{->}(0.0,1.0)(5.0,0.0)
\psline[linewidth=0.2mm]{->}(0.0,4.0)(4.0,5.0)
\psline[linewidth=0.2mm]{->}(4.0,5.0)(5.0,0.0)
\psline[showpoints=true,linewidth=0.6mm]{-}(2.25,2.5)(2.25,2.5) 
\psline[linewidth=0.4mm]{->}(4.5,2.5)(0,2.5)
\psline[linewidth=0.4mm]{->}(2.0,4.5)(2.5,0.5)
\psline[showpoints=true,
linewidth=0.4mm]{->}(4.5,2.5)(8.751,3.400) 
\psline[showpoints=true,linewidth=0.4mm]{->}(0.0,2.5)(-1.5,2.5)
\psline[showpoints=true,linewidth=0.4mm]{->}(2.0,4.5)(1.287,7.352)
\psline[showpoints=true,linewidth=0.4mm]{->}(2.5,0.5)(1.650,-3.751)
\rput(1.0,3.0){$_{_{0}}b$}
\rput(3.0,1.5){$_{_{1}}b$}
\rput(6.6,3.6){$_{_{0}}f$}
\rput(-1.2,3.2){$_{_{1}}f$}
\rput(2.2,6.4){$_{_{2}}f$}
\rput(2.7,-1.4){$_{_{3}}f$}
\rput(5.4,2.0){$_{_{0}}Z^{\,p}$}
\rput(-0.8,1.9){$_{_{1}}Z^{\,p}$}
\rput(1.0,4.7){$_{_{2}}Z^{\,p}$}
\rput(1.6,0.0){$_{_{3}}Z^{\,p}$}
\rput(2.9,3.0){$Z^{\,n}$}
\end{pspicture}
\caption{\textsl{Face vectors and port locations
(nodal section).\hfill}}\label{F:3}
\end{figure}
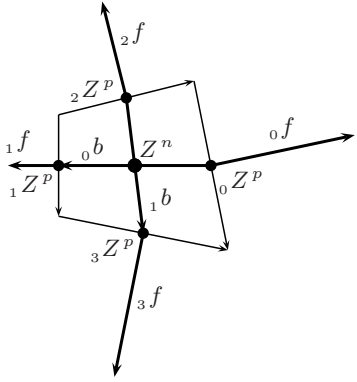

At every cell face ${\iota\in\{0,...,5\}}\/$ and for any 
given $\tau\in\mathbb{R}_{+}\,$ the following time shifted finite
differences of $\,Z\,$ in directions ${ _{\mu} b }$ (\,$\mu = 0,1,2\,$)
form a vector valued function
\begin{equation}\label{3.2}\centering
\begin{split}
_{\iota}\!{\nabla}^{B} Z_{\mu}\,(\,t\,)\;:\,=\;
\begin{cases}
\,2\,(-1)^{\iota}(\,Z^{n}\, _{\mid\,t-\tau/2} 
-\,_{\iota}Z^{p}\,_{\mid\,t}\,)\quad
&\text{if $\mu\,=\,[\iota / 2]$}\\
\,(\,\,_{2\mu+1} Z^{p}\,-\,_{2\mu} Z^{p}\,\,)
\,_{\mid\,t -\tau\,}\quad
&\text{if $\mu\,\neq\,[\iota / 2]$}\,
\end{cases}
\end{split}
\end{equation}
($\,[\,x\,]$ denotes the \emph{integer part} of $\,x\in\mathbb{R}\,$).
The time increments are chosen conform with the updating conventions
of DSC schemes (as will be seen in a moment) and are consistent.
In fact, in the first order of the time increment ${\,\tau\,}$
and of the linear cell extension,
the vector ${\,_{\iota}\!{\nabla}^{B} Z\,}$
in the centre point of face ${\,\iota\,}$ approximates
the scalar products of the node vectors with the gradient
${\,\nabla Z\,}$.
Let, precisely, for a fixed centre point on face ${\,\iota\,}$
and $\,\epsilon \in \mathbb{R}_{+}\,$ the \emph{$\epsilon$-scaled cell}
have edge vectors
$\,_{\iota}e\sptilde\,:\,=\,\epsilon\,\,_{\iota}e\,$. 
Let also $\,_{\iota}{\nabla}^{B\sptilde}Z_{\mu}\,$ denote
function \eqref{3.2} for the $\epsilon$-scaled cell (with node vectors 
$\,_{\mu}b\sptilde\,=\,\epsilon\,_{\mu}b\;$).
Then at the fixed point holds
\begin{equation}\label{3.3}\centering
\begin{split}
<\,_{\mu}b\,,\,\text{grad($Z$)}\,>\,\,\,
=\,\,_{\mu}b\cdot\nabla Z\,\,
=\,\,\lim_{\epsilon\to 0}\,\,\lim_{\tau \to 0} \, \,
\frac{1}{\epsilon}\,_{\iota}\!{\nabla}^{B^{\sptilde}}Z_{\mu} \, ,
\end{split}
\end{equation}
as immediately follows from the required $C^1$-smoothness of the field $Z$.

To recover the gradient
${{\nabla}Z\,}$ from \eqref{3.2}
in the same order of approximation,
observe that for every orthonormal
basis $\,{(_{\nu}u)_{\nu=0,...,m-1}}\,$ of
$\mathbb{R}^{m}\,\text{or}\,\,\mathbb{C}^{m}\,$, and for
any basis $\,{(_{\mu}b)_{\mu = 0,...,m-1}}\,$ with coordinate
matrix ${\beta_{\nu}^{\mu}}\,=\,{<\, _{\nu}u\,,\, _{\mu}b\,>}$,
the scalar products of every vector $\,a\,$ with $\,{_{\mu} b}\,$ equal
\begin{equation}\label{3.4}\centering
\underbrace{<\,_{\mu}b\,,\,a\,>}_{\qquad=\,: 
\,\,{\alpha}_{\mu}^{B}}\,\,=\,\sum\nolimits_{\nu=0}^{m-1}\,
\underbrace{<\,_{\mu}b\,,\,_{\nu}u\,>}_{\,\,\,
({\bar{\beta}}_{\mu}^{\nu})\,=\,({\beta}_{\nu}^{\mu} )^{^{*}}}\,
\underbrace{<\,_{\nu}u\,,\,a\,>}_{\qquad=\,:\,\,{\alpha}_{\nu}}
\,\,=\,\bar{\beta}_{\mu}^{\nu}\,{\alpha}_{\nu}\;
\end{equation}
(\,at the right-hand side -~and henceforth~- observe
\textsc{Einstein}'s convention\,),
hence
\begin{equation}\label{3.5}\centering
{\alpha}_{\nu}\,=\,
{\gamma}_{\nu}^{\mu}\alpha_{\mu}^{B}\,,
\qquad\text{with}\qquad({\gamma}_{\nu}^{\mu}) 
\,:\,=\,{({(\beta_{\nu}^{\mu})}^{*})}^{-1}\quad .
\end{equation}
In words: The scalar products of any vector with the basis
vectors ${_{\mu} b\,}$ transform into the coordinates of that vector
with respect to an orthonormal basis $\,{_{\nu}u}\,$ by multiplication
with matrix $\,{\gamma=(\beta^*)^{-1}}\,$, where
$\,{\beta_{\nu}^{\mu}}\,=\,{<\,_{\nu}u\,,\,_{\mu}b\,>}\,$,\,
i.e. $\,\beta\,$ is the matrix of the coordinate (column) vectors
$\,{_{\mu}b}\,$ with respect to the given ON-basis
$\,{_{\nu}u}\,$,\, and $\,\gamma\,$ its adjoint inverse.

This applied to the node vector basis ${ _{\mu}b\,}$ and \eqref{3.3}
yields the approximate gradient of $\,Z\,$ at face $\iota$
\begin{equation}\label{3.6}\centering
_{\iota}\!\nabla Z_{\nu}\quad
=\quad{\gamma}_{\nu}^{\mu}\,\,\,_{\iota}\!{\nabla}^{B}Z_{\mu}.
\end{equation}
The scalar product of the gradient with face vector  
${_{\iota}f^{\nu}}\,=\,{<\,_{\iota}f,\,_{\nu}u>}\,$,\\
$\,\nu\in\{0,1,2\}\,$ is thus
\begin{equation}\label{3.7}\centering
\begin{split}
\begin{aligned}
\qquad _{\iota}S\; 
&=\;_{\iota}f\,\cdot\,_{\iota}\!\nabla Z\;
=\underbrace{_{\iota}f^{\nu}\;
{\gamma}_{\nu}^{\mu}}_{\qquad\;=\,:\,\,_{\iota} s^{\mu}}\,
_{\iota}\!{\nabla}^{B} Z_{\mu}\,\,
=\,\,_{\iota} s^{\mu}\,\,_{\iota}\!{\nabla}^{B} Z_{\mu}\;.
\end{aligned}
\end{split}
\end{equation}
Continuity of the gradient at cell interfaces yields
linear updating equations for $Z^p$ on the two adjacent faces.  
In fact, for any two neighbouring cells
$\zeta$, $\chi$ with common face, labelled $\iota$
in cell $\zeta$ and $\kappa$ in $\chi\,$, continuity requires
\begin{equation}\label{3.8}\centering
_{\iota}^{^{\zeta}}\!S\quad=\quad-\,\,_{\kappa}^{^{\chi}}\!S\,.
\end{equation}
Substituting \eqref{3.7} for
$\,_{\iota}^{^{\zeta}}\!S\,$ and $\,_{\kappa}^{^{\chi}}\!S\,$
and observing the time shifts in \eqref{3.2}
\linebreak
provides the updating relations for 
$\,Z^{\,p}\,$ at the cell interfaces.\\
To derive these explicitely, we first introduce the following quantities
$\,{_{\iota}z_{\mu}^{\,p,\,n}}\,$, ($\,\iota\,=\,0,...,5\,$;
$\mu\,=\,0,1,2\,$)
\begin{equation}\label{3.10}\centering
\begin{split}
_{\iota}z_{\mu}^{n}\,(\,t\,)\quad :\,=\quad
\begin{cases}
\,\,2\,(-1)^{\iota}\,\,Z^{\,n}\,_{\mid\,t}\qquad
&\text{if $\mu\,=\,[\iota /2]$}\,\,\\
\,\,(\,_{2\mu +1} Z^{\,p} 
-\,_{2\mu} Z^{\,p}\,)_{\mid\,t-\tau/2}\qquad
&\text{else}\,
\end{cases}\;,
\end{split}
\end{equation}
which in virtue of \eqref{1.1} yields
$\;{_{\iota}z_{\mu}^{\,p}}\,=\,{(\,p\,,\,Z\,)}\,
=\,{(\,p\sptilde,\,Z\circ\,_{\iota}{\sigma}^{-1}\,)}\,
=\\=\,{_{\iota}z_{\mu}^{\,n}\,\mid{Z\,\circ\,_{\iota}{\sigma}^{-1}}}\,$,
where $\,_{\iota}{\sigma}\,:\,n\,\mapsto\,p\,$
denotes the nodal shift pertinent to face $\,\iota\,$.\;
In particular
\begin{equation}\label{3.11}
_{\iota}z_{[\iota/2]}^{\,p}\,(\,t\,)\quad=\quad
\,\,2\,(-1)^{\iota}\,\,_{\iota}Z^{\,p}\,_{\mid\,t}\;,
\end{equation}
which together with \eqref{3.10} for $\,\mu\,\neq\,[\iota/2]\,$
is consistent with
\begin{equation}\label{3.12}
_{\iota}z_{\mu}^{\,n}\,(\,t+\tau/2\,)\quad
=\quad -\;\frac{1}{2}\,(\,_{\,2\mu+1}z_{\mu}^{\,p}\,
+\, _{\,2\mu}z_{\mu}^{\,p}\,)\,(\,t\,)\;.
\end{equation} \newline
From (\,\ref{3.2}, \ref{3.7}, \ref{3.10}, \ref{3.11}\,) follows that
\begin{equation}\label{3.13}
\begin{split}
_{\iota}S\,_{\mid\,t+\tau}\quad
&=\quad\,_{\iota}s^{\mu}\,(\,_{\iota}z_{\mu}^{\,n}\,_{\mid\,t+\tau/2}\,
-\,2\,{(-1)}^{\iota} {\delta}_{\mu}^{[\iota/2]}\,\,
_{\iota}Z^{\,p}\,_{\mid\,t+\tau}\,) \\
&=\quad _{\iota}s^{\mu}\,(\,_{\iota}z_{\mu}^{n}\,_{\mid\,t+\tau/2} \,
-\,{\delta}_{\mu}^{[\iota/2]}\,\, 
_{\iota}z_{\mu}^{\,p}\,_{\mid\,t+\tau}\,)
\quad .
\end{split}
\end{equation}
This, with (\ref{3.8},\ref{3.10})
and the continuity of $\,Z\,$, i.e.
$\,_{\iota}^{^{\zeta}}Z\,^{p}\,=\,_{\kappa}^{^{\chi}}Z\,^{p}\,$, implies
\begin{equation}\label{3.14}
_{\iota}^{^{\zeta}}z\,_{[\iota/2]}^{p}\,(\,t+\tau\,)\,=\,
\,\frac{\,_{\iota}^{^{\zeta}}s\,^{\mu}\,\, 
_{\iota}^{^{\zeta}}z\,_{\mu}^{n}\,(\,t+\tau/2\,)\,
+\,_{\kappa}^{^{\chi}}s\,^{\nu}\,\,\,
_{\kappa}^{^{\chi}}z\,_{\nu}^{n}\,(\,t+\tau/2\,)} 
{_{\iota}^{^{\zeta}}s\,^{[\iota/2]}
+\,( -1 )^{\iota +\kappa}\;_{\kappa}^{^{\chi}}s\,^{[\kappa/2]}}\;,
\end{equation}
and for completeness we agree upon setting
$\,_{\iota}^{^{\zeta}}z\,_{\mu}^{p}\,(\,t+\tau\,)\,
:\,=\,_{\iota}^{^{\zeta}}z\,_{\mu}^{n}\,(\,t+\tau/2\,)\,$,
for $\,\mu\,\neq\,[\iota/2]\,$.
Note that the latter relations contain a slight inconsistency,
in that continuity might be infringed -~which is yet circumvented
by taking the arithmetic means of the two adjacent values.
In fact, our agreement doesn't do harm, since any discontinuity
disappears with mesh refinement. 

We have thus defined complete recurrence relations for
$\,{z^{\,p}}\,$ (\,given $\,{z^{\,n}}\,$ on the former reflection step),
which at the same time determine on face ~$\iota$
the field components and their gradients
\begin{equation}\label{3.15}\centering
_{\iota}\!\nabla Z_{\nu}\quad
=\quad{\gamma}_{\nu}^{\mu}\; _{\iota}z_{\mu}^{\,p},
\end{equation}
and which essentially constitute the connection step of the algorithm.
\newline
Nodal gradients are similarly, yet more simply, derived using
\begin{equation}\label{3.16}\centering\notag
{\nabla}^{B} Z_{\mu}^{\,n}\,(\,t\,+\frac{\tau}{2}\,)\quad:\,
=\quad(\,_{2\mu+1} Z^{\,p}\,-\,_{2\mu} Z^{\,p}\,)\,(\,t\,)\;;
\quad\mu\,=\,0,\,1,\,2\;
\end{equation}
in the place of \eqref{3.2} and then again \eqref{3.6}.
With the node and cell-boundary values and gradients of $\,T\,$ and
$\,\vec{u}\,$, the nodal updating relations for these quantities
are directly derived from equations~(\,\ref{2.4}, \ref{2.5}\,)
in sect.~\ref{S:sec2}\,. For equations \eqref{2.4} (without the
convective term) this has essentially been carried out in \cite{He1},
sect.~5, and the procedure remains straightforward.

The obtained updating relations are explicit and consistent
with near-field interaction (\,only adjacent quantities enter\,).
So, they can optionally be transformed into scattering relations
for incident and reflected quantities \eqref{1.6} along the
guidelines of section~\ref{S:sec1} -~with established advantages
for the stability estimates \cite{He2}.
\vspace{-0.2cm}
\section{Pressure}\label{S:sec4}
Pressure is a known subject sui generis in Computational Fluid
Dynamics \cite{ATP}\cite{MeSt}\cite{GDN}
already insofar as pressure fluctuations typically do not match
the time scales of heat propagation and fluid flow.
Pressure fluctuations are related to acoustic waves,
the net effect of which can be important and usually has a strong
impact on computational stability \cite{LMDM}.

Pressure fluctuations play also a key role in the following procedure,
which is known as \emph{divergence cleaning} in Magnetohydrodynamics
[ibid., p.~128] and ensures mass conservation in the present context.

For Boussinesq-incompressible fluids, conservation of mass simply requires
divergence-free flow, $\;{div\,\vec{u}\,=\,0}\;$,
i.e. in integral form, using Gauss' theorem,
$0\,=\,\int\nolimits_{\,\zeta}\,div\;\vec{u}\;dV\,
=\,\int\nolimits_{\,\partial\zeta}\,\vec{u}\,\cdot\,dF\,$.
Since equations (\,\ref{2.1}, \ref{2.2}\,) in section~\ref{S:sec2}
do not a priori ensure this, additional arrangements must be
made -~which is done by means of pressure.

In a successive overrelaxation (SOR) routine, interposed between
the connection and reflection steps of the iteration cycle,
firstly the (discrete) right-hand side integrals
$\,I_{\partial\zeta}\,=\,\int _{\partial\zeta}\vec{u}\cdot dF\,$
are computed, and then the pressure $\,p\,$
which compensates $\,I_{\partial\zeta}\,$ so that
\begin{equation}\centering\label{4.2}
\quad\frac{\tau}{\varrho_{\infty}}\,
\int\nolimits_{\,\partial\zeta}\, grad\,p\,\cdot\,dF\quad
=\quad\int\nolimits_{\,\partial\zeta}\,\vec{u}\,\cdot\,dF\;.
\end{equation}
Indeed, taking $\,Z\,$ of the preceeding section as the pressure,
it follows from (\,\ref{3.10}, \ref{3.12}\,) that equations \eqref{4.2}
(\,in discrete form with sums over the cell faces, of course\,)
yield a unique solution $\,p^{\,n}\,=\,Z^{\,n}\,$ for every cell,
given the right-hand side integral $\,I_{\partial\zeta}\,$.
Note that we are actually solving Poisson's equation
$\;{\Delta p}\,=\,{(\,\varrho_{\infty}/\tau)}\,{div\,\vec{u}}\;$
in integral form.
\newline
With the new face pressure gradient,
extracted from (\,\ref{3.10},\,\ref{3.14},\,\ref{3.15}\,),
the face values of $\,\vec{u}\,$ are updated as
$\,\vec{u}\,-\,(\tau/\varrho_{\infty})\,grad\,p\,$.

After each SOR cycle, continuity of $\,grad\;p\,$
at the cell faces must be re-established using \eqref{3.8},
i.e. by updating the port values of $\,p\,$ according to \eqref{3.13}.
The chain of processes is reiterated until
$\sum \nolimits_{\zeta}\,I_{\partial\zeta}\,<\,\epsilon\,$
for a suitable bound $\,\epsilon\,$
(\,which happens after a few iterations for appropriate choices\,).
\vspace{-0.2cm}
\section{Coaxial line}\label{S:sec5}
To illustrate the approach in a stalwart application, 
we display the results of simulations with coaxial line RL230-100
under high power operating conditions
(\,realistically inferred from a ion cyclotron resonance heating 
\emph{ICRH} experiment in plasma physics\,).

The inner and outer conductors of diameters 100 mm and 230 mm
are made of copper and aluminium, respectively, and the rigid
line is filled with air at atmospheric pressure.
We have simulated the heating process from standby to steady state
CW operation, at frequency 100 MHz and 160 kW transmitted power,
for horizontal position of the line and with outer conductor cooled
at 40 degrees Celsius.

\begin{figure}[!h]\centering
\setlength{\unitlength}{1.cm}
\begin{pspicture}(0.0,0.0)(8.0,5.8)\centering
\psset{xunit=1.0cm,yunit=1.0cm}
\rput(0.50,3.70){\includegraphics[scale=0.80,clip=0]{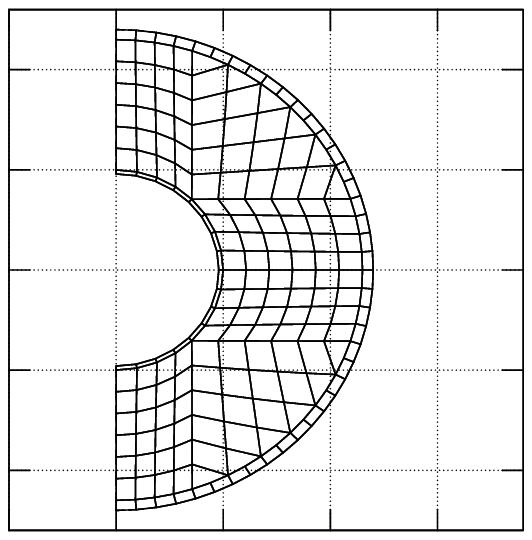}}
\rput(-0.90,1.25){\small (a)}
\rput(6.60,3.70){\includegraphics[scale=0.80,clip=0]{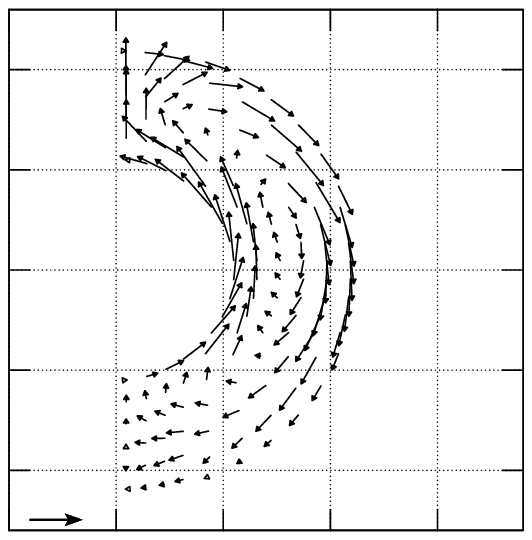}}
\rput(5.25,1.25){\small (b)}
\end{pspicture}
\vspace{-0.8cm}
\caption{
\textsl{
(a) mesh 
(b) velocity profile [ reference arrow: 0.1 $m s^{-1}\,$ ]
\hfill}}
\label{F:4}
\end{figure}

Figure~\ref{F:4}~b displays the computed air flow profile (vertical section)
in steady state, which is attained some minutes after switch-on.
Visibly, the natural convection pattern is nicely developed.

The computations have been carried through with a 3D-mesh of 10 layers
in axial direction, over 200 millimeters of line, the transverse cross
section of which is displayed in figure~\ref{F:4}~a.
At the metallic interfaces no-slip boundary conditions are
implemented and free-slip conditions at all other boundaries.

Simultaneously, a Maxwell field TLM algorithm was run to provide
the heat sources.

\vspace{-0.2cm}
\section{Conclusions}\label{C:sec6}
A prototypical implementation of the \textsc{Oberbeck-Boussinesq}
approximation to viscous flow has been presented in this paper,
which demonstrates the fundamental fitness of the DSC approach for
fluid dynamic computations.
In this respect, at least, (this was recently called into question
by a TLM expert) DSC schemes significantly transcend the range of
application of the TLM method, from which they descend.

A next natural step in the line of this study is, of course, the
implementation of turbulence models which are compatible with the
Boussinesq approach, such as the $\,k-\epsilon\,$ model \cite{ATP}\,,
first of all. -~We hope this paper stimulates some interest into
joint further investigation in that direction.
\vspace{-0.2cm}
\vspace{25pt}\
\textbf{Acknowledgement}
\newline
The author wants to thank Ingolf Lehniger for some stylistic improvements.
\vspace{-0.2cm}

\textsc{Spinner} GmbH. M\"unchen; Aiblinger Str. 30, DE-83620 Westerham
\newline
E-mail address:\; s.hein@spinner.de 

\vspace{5pt}
\begin{thebibliography}{99}

\bibitem[MeSt]{MeSt}
Meister, A., Struckmeier, J.,
{\it Hyperbolic Partial Differential Equations\/}
Theory, Numerics and Applications,
Friedrich Vieweg and Sohn, G\"ottingen 2002

\bibitem[LeVeque]{LMDM}
LeVeque, R.J., Mihalas, D., Dorfi, E.A., M\"uller, E.
{\it Computational Methods for Astrophysical Fluid Flow\/},
Saas Fee Advanced Courses, 27,
Springer-Verlag Berlin Heidelberg, 1998

\bibitem[GDN]{GDN}
Griebel, M., Dornseifer, T., Neunhoeffer, T., {\it Numerical Simulation in
Fluid Dynamics\/}, SIAM monographs on mathematical modeling and computation,
Society for Industrial and Applied Mathematics, 1998

\bibitem[ATP]{ATP}
Anderson, D.A., Tannehill, J.C., Pletcher, R.H.,
{\it Computational Fluid Mechanics and Heat Transfer\/},
series in computaional methods in mechanics and thermal sciences,
Hemisphere Publishing Corporation, 1984

\bibitem[Bss]{Bss}
Boussinesq, J.,
{\it Th\'eorie Analytique de la Chaleur\/},
Gauthiers-Villars, 2., Paris 1903

\bibitem[Obb]{Obb}
Oberbeck, A.,
\"Uber die W\"armeleitung der Fl\"ussigkeiten bei Ber\"ucksichtigung
der Str\"omung infolge Temperaturdifferenzen.,
Ann. Phys. Chem., vol. 7, pp. 271-292, 1879

\bibitem[JoB]{JoB}
Johns, P.B., Beurle R.L., Numerical solution of 2-dimensional scattering
problems using transmission line matrix, Proc. IEEE, vol. 118,
pp. 1203-1208, 1971

\bibitem[Tlm]{Tlm}
Proceedings of the 2nd Int. Workshop on Transmission Line Matrix (TLM)
Modelling, TU M\"unchen, 1997

\bibitem[He0]{He0}
Hein, S., Finite-difference time-domain approximation of Maxwell's equations
with nonorthogonal condensed TLM mesh, Int. J. Num. Modelling, 
vol. 7, pp. 179-188, 1994 

\bibitem[He1]{He1}
\makebox{Hein, S., Dual scattering channel schemes extending the
\textsc{Johns} Algorithm}
\newline http://arxiv.org/abs/math.NA/0309261, March 2004

\bibitem[He2]{He2}
\makebox{Hein, S., On the stability of dual scattering channel schemes,}
\newline http://arxiv.org/abs/math.NA/0405095, preprint, May 2004

\bibitem[He3]{He3}
Hein, S., TLM numerical solution of Bloch's equations for magnetized gyrotropic
media, Appl. Math. Modelling, vol. 21, pp. 221-229, 1997

\end{thebibliography}
\end{document}